\documentclass[12pt,reqno,a4wide]{amsart}

\usepackage{tikz} 
\usepackage{calrsfs}
\usepackage[charter]{mathdesign}

\allowdisplaybreaks

\title{Summing a family of generalized Pell numbers}

\author{Helmut Prodinger}
\address{Helmut Prodinger\\
	Department of Mathematical Sciences\\
	Stellenbosch University\\
	7602 Stellenbosch\\
	South Africa}
\email{hproding@sun.ac.za}


\begin{document}
\begin{abstract}
	A new family of generalized Pell numbers was recently introduced and studied by  Br\'od \cite{Dorota}. These number possess, as Fibonacci numbers, a Binet formula. Using this, partial sums of arbitrary powers of generalized Pell numbers can be summed explicitly. For this, as a first step, a power $P_n^l$ is expressed as a linear combination of $P_{mn}$.
	The summation of such expressions is then manageable using generating functions. Since the new family contains a parameter $R=2^r$, the relevant manipulations are quite involved, and computer algebra produced huge expressions that where not trivial to handle at times.
		\end{abstract}

\maketitle
	
	\section{Introduction}

	Br\'od \cite{Dorota} recently introduced $r$-Pell numbers, via
	\begin{equation*}
P(r,n)=C_1r_1^n+C_2r_2^n,
	\end{equation*}
	with
	\begin{align*}
r_{1,2}=\frac{2^r\pm\sqrt{4^r+2^{r+1}}}{2},\quad C_{1,2}=1\pm\frac{2^r+1}{\sqrt{4^r+2^{r+1}}}.
	\end{align*}
	The classical Pell numbers are then given as $P_{n+2}=P(1,n)$. First, we write $R=2^r$, and $W=\sqrt{R^2+2R}$.
	Then we have
	\begin{equation*}
C_1r_1^n+C_2r_2^n=\frac{2}{RW}(r_1^{n+2}-r_2^{n+2}),
	\end{equation*}
	and we will consider the shifted versions
	\begin{equation*}
P_R(n):=\frac{2}{RW}(\lambda^{n}-\mu^{n}), \quad\text{ with }\quad \lambda=\frac{R+\sqrt{R^2+4R}}{2},\ \mu=\frac{R-\sqrt{R^2+4R}}{2}.
	\end{equation*}
	This version is slightly more appealing, and fits better with the classical Pell numbers: $P_2(n)=P(n)$.
	(For $R=1$ we are in the world of Fibonacci numbers.) The recursion is 
\begin{equation*}
P_R(n+2)-2RP_R(n+1)-RP_R(n)=0,
\end{equation*}
and the first few numbers are
\begin{equation*}
0,\frac2{R} ,2,2R+1,2R \left( R+1 \right) ,\frac R2 \left( 4{R}^
{2}+6R+1 \right) ,\frac{{R}^{2}}2 \left( 2R+3 \right)  \left( 2R+1
\right),\dots
\end{equation*}

We also need the associated sequence
\begin{equation*}
Q_R(n):=\frac{2}{R}(\lambda^{n}+\mu^{n}),
\end{equation*}
and it is easy to see that $Q_R(n)=2P_R(n+1)-RP_R(n)$.

	The goal of the present note is to compute the sums
	\begin{equation*}
\sum_{0\le k \le n}P_R(km)^\ell
	\end{equation*}
	efficiently, for positive parameters $m$ and $\ell$.
	
	We presented such a treatment in \cite{Prodinger-Fro} for Fibonacci numbers, and, more recently, also in \cite{Prodinger-balancing} for balancing numbers studied by Komatsu \cite{Takao}, but the current instance is definitely the most challenging, as it contains the parameter $R$, and that makes computations difficult (even for a computer).
	
	\section{Generating functions}
	
	First we give the generating function for the generalized Pell numbers, which follows immediately from the Binet formula:
	\begin{equation*}
\sum_{n\ge0}P_R(n)z^n=\frac{4z}{R(2-2zR-z^2R)}.
	\end{equation*}
	We have
	\begin{multline*}
	2^mP_R((n+2)m)
	-(R(R+1)2^{m-1}P_R(m)- R^22^{m-2}P_R(m-2))P_R((n+1)m)\\*
	+(-1)^mR^mP_R(nm)	=0.
	\end{multline*}
For $m=1$, this is just the ordinary recursion for generalized Pell numbers, but we need this more general version.

	This formula is not too difficult to prove using the Binet formula and the help of a computer, but not so easy
	to find. We used Maple's package \textsf{gfun} for that.
	Let us rewrite this in terms of generating functions: Set
	\begin{equation*}
F_m(z)=\sum_{n\ge0}P_R(mn)z^n,
	\end{equation*}
	then
	\begin{multline*}
2^m\sum_{n\ge2}P_R(nm)z^n-z\Bigl(R(R+1)2^{m-1}P_R(m)-R^22^{m-2}P_R(m-2)\Bigr)\sum_{n\ge1} P_R(nm)z^m
	\\*+z^2(-1)^mR^m\sum_{n\ge0}P_R(nm)z^n=0
	\end{multline*}
	or
	\begin{multline*}
	2^m(F_m(z)-zP_R(m))-z \Bigl(R(R+1)2^{m-1}P_R(m)-R^22^{m-2}P_R(m-2)\Bigr)F_m(z) 
	\\*+z^2(-1)^mR^m F_m(z)=0,
	\end{multline*}
	or
	\begin{equation*}
F_m(z)=\frac{ 2^mzP_R(m)}{ 2^m-z \bigl(R(R+1)2^{m-1}P_R(m)-R^22^{m-2}P_R(m-2)\bigr) 
	+(-1)^mz^2R^m }.
	\end{equation*}
	For $m=1$, we get the first generating function again as a check.
	
	Now let us abbreviate
	\begin{equation*}
	F_m(z)=\frac{ az}{ b+cz+dz^2},
	\end{equation*}
	then
	\begin{equation*}
	\frac1{1-z}F_m(z)=\frac{ a}{ b+c+d}\frac1{1-z}-\frac{ a}{ b+c+d}\frac{ b-dz}{ b+cz+dz^2},
	\end{equation*}
	and, by reading off coefficients,
	\begin{align*}
[z^n]\frac1{1-z}F_m(z)=\sum_{k=0}^{n}P_R(km)&=\frac{ a}{ b+c+d}\\
&+\frac{ 1}{ b+c+d}\Bigl(bP_R(m(n+1))+dP_R(mn)\Bigr).
	\end{align*}

	\section{Conversion of powers}
We need to rewrite $P_R(n)^m$ in a `linearized' form. The cases $m$ odd resp.\ even are significantly different in complexity.
\begin{equation*}
P_R(n)^{2m+1}=\frac{2^{2m}}{R^{3m}(R+2)^{m}}\sum_{j=0}^m(-1)^j\binom{m}{j}
P_R((m-2j)n)\Bigl(-\frac R2\Bigr)^{jn}
\end{equation*}
The first one is basically a rearrangement following the binomial theorem:
\begin{align*}
(x-x^{-1})^{2m+1}&=\sum_{j=0}^{2m+1}\binom{2m+1}{j}(-1)^jx^{2j-2m-1}\\
&=\sum_{j=0}^{m}\binom{2m+1}{j}(-1)^jx^{2j-2m-1}
+\sum_{j=m+1}^{2m+1}\binom{2m+1}{j}(-1)^jx^{2j-2m-1}\\
&=\sum_{j=0}^{m}\binom{2m+1}{j}(-1)^jx^{2j-2m-1}
-\sum_{j=0}^{m}\binom{2m+1}{j}(-1)^jx^{-2j+2m+1}\\
&=\sum_{j=0}^{m}\binom{2m+1}{j}(-1)^j(x^{2j-2m-1}-x^{-2j+2m+1}).
\end{align*} 
 The second one does not follow such a simple symmetry:
 \begin{align*}
P_R(n)^{2m}&=\frac{2^{2m+1}}{R^{3m}(R+2)^m}\sum_{j=0}^{m-1}(-1)^j\binom{2m}{j}
\frac{P_R(2(m-j)(n+1))}{P_R(2(m-j)}\Bigl(-\frac R2\Bigr)^{jn}\\
&-\frac{2^{2m-1}}{R^{3m-1}(R+2)^m}\sum_{j=0}^{m-1}(-1)^j\binom{2m}{j}
\frac{P_R(2(m-j)n)Q_R(2(m-j))}{P_R(2(m-j)}\Bigl(-\frac R2\Bigr)^{jn}\\
&+(-1)^m\binom{2m}{m}\frac{2^{2m}}{R^{3m}(R+2)^m}\Bigl(-\frac R2\Bigr)^{mn}.
\end{align*}
This identity was found by creative guessing with a computer, taking two days. Using the Binet forms, this is basically a polynomial
identity, like the one that we used earlier in \cite{Prodinger-variants}:
\begin{equation*}
x^{2N}-x^{-2N}=(x+x^{-1})\sum_{l=0}^{N-1}\binom{N+l}{N-l-1}(x-x^{-1})^{2l+1}.
\end{equation*}
We invite the readers to provide a formal proof, which should be basically and exercise in manipulating binomial coefficients.
The creative part of this was to \emph{find} the identity.

 \section{Summing powers}
 
 Using the formulae of the previous section, we can evaluate
 \begin{equation*}
\sum_{0\le k \le n}P_R(\ell k)^{m}.
 \end{equation*}
 We have shown already how to sum
 \begin{equation*}
 \sum_{k=0}^{n}P_R(km).
 \end{equation*}
 We need the slightly more general form, with a parameter $\sigma$,
 \begin{equation*}
 \sum_{k=0}^{n}P_R(km)\sigma^k,
 \end{equation*}
 because of the extra powers $(-R/2)^{jn}$ in our formulae. But that is not difficult, since
 in 
 \begin{equation*}
 F_m(z)=\frac{ 2^mzP_R(m)}{ 2^m-z \bigl(R(R+1)2^{m-1}P_R(m)-R^22^{m-2}P_R(m-2)\bigr) 
 	+(-1)^mz^2R^m }=\frac{ az}{ b+cz+dz^2}
 \end{equation*}
 we just have to replace $z\to \sigma z$, and
 \begin{equation*}
 \frac1{1-z}F_m(\sigma z)=\frac{ a\sigma}{ b+c\sigma+d\sigma^2}\frac1{1-z}-\frac{ a\sigma}{ b+c\sigma+d\sigma^2}\frac{ b-d\sigma^2z}{ b+c\sigma z+d\sigma^2z^2},
 \end{equation*}
 and, reading off coefficients again,
 \begin{align*}
 [z^n]\frac1{1-z}F_m(\sigma z)=\frac{ a\sigma}{ b+c\sigma+d\sigma^2}
 -\frac{1}{ b+c\sigma+d\sigma^2}\bigl(-d\sigma^2T(\sigma n)+bT(\sigma (n+1))\bigr).
 \end{align*}
 
 \section{Examples}
 
 Here is the generating function of the numbers
 \begin{equation*}
\sum_{k=0}^nP_R(k)^2:
 \end{equation*}
 \begin{multline*}
16 {\frac {R-2}{ \left( 1-z \right) {R}^{2} \left( R+2 \right) ^{2}
		\left( 3 R-2 \right) }}
+16 {\frac {-zR+z{R}^{2}+2}{ \left( R+2 \right)  \left( 3 R-2
		\right) {R}^{2} \left( {z}^{2}{R}^{2} -4 zR(R+1)+4 \right) }}\\*
-16 {\frac {1}{ \left( zR+2 \right) {R}^{2} \left( R+2 \right) ^{2}}}
 \end{multline*}
 The coefficients can be expressed with generalized Pell numbers, as described in the previous sections.
 
 Another example is this:
 \begin{equation*}
 \sum_{k=0}^nP_R(2k)^3,
 \end{equation*}
 with generating function
 \begin{multline*}
{\frac {512({R}^{6}+32 {R}^{4}+32 {R}^{3}+64)}{ \left( R+2 \right) ^
		{2} \left( 3 R-2 \right)  \left( 7 {R}^{2}-2 R+4 \right)  \left( 3
		 {R}^{2}+6 R+4 \right)  \left( {R}^{2}+2 R-4 \right)  \left( {R}^{3
		}-4 {R}^{2}-8 \right)  \left( z-1 \right) }}
\\-{\frac {384(z{R}^{6}-64)}{ \left( {z}^{2}{R}^{6}-16 z{R}^{3}(R+1)  +64\right) R \left( {R}^{3}-4 {R}^{2}-8 \right)  \left( {R}^{2}
		+2 R-4 \right)  \left( R+2 \right) ^{2}}}
\\-{\frac {128 \left( 2 R+3 \right)  \left( 2 R+1 \right)  \left( z{
			R}^{6}-64 \right) }{ \left( {z}^{2}{R}^{6}-16zR^3(R+1)(4R^2+8R+1)+64 \right) R \left( 3 {R}^{2}+6 R+4
		\right)  \left( 7 {R}^{2}-2 R+4 \right)  \left( 3 R-2 \right) 
		\left( R+2 \right) ^{2}}}
 \end{multline*}
 
 \bigskip
 We hope that this paper will trigger more research related to this new exciting families of numbers.
 
	\bibliographystyle{plain}

\end{document}